\documentclass{elsart3-1}

\usepackage{mathtools}                                                     
\mathtoolsset{showonlyrefs=true}                                    

\usepackage{amssymb}
\usepackage{graphicx}
\usepackage{amsmath}
\usepackage{color}

\usepackage[english,francais]{babel}

\def \a {{\alpha}}

\def \d {{\delta}}

\def \R {{\mathbb {R}}}
\def \N {{\mathbb {N}}}

\def \Rdd{{\mathbb{R}\times\mathbb{R}^{d}}}
\def \Rpz {{\mathbb {R}^{p_0}}}
\def \Rpo {{\mathbb {R}^{p_1}}}

\def \Rpp{{\mathbb{R}\times\mathbb{R}^{p_0}\times\mathbb{R}^{p_1}}}

\def \x {{\xi}}

\def \eps {{\varepsilon}}

\def \y {{\eta}}

\def \z {{\zeta}}

\def \tt {{\bar t}}
\def \g {{\gamma}}
\def \O {{\Omega}}
\def \phi {{\varphi}}

\def\p{\partial}

\def \è {\`e }

\newcommand\Lc{\mathcal{L}}
\newcommand\Gc{\mathcal{G}}

\newcommand\Nb{\mathbb{N}}
\newcommand\Nzero{\mathbb{N}_{0}}
\newcommand{\norm}[1]{\| #1 \|}

\newtheorem{theorem}{Theorem}[section]
\newtheorem{lemma}[theorem]{Lemma}
\newtheorem{e-proposition}[theorem]{Proposition}

\newtheorem{e-definition}[theorem]{Definition\rm}

\newtheorem{notation}[theorem]{Notation}
\newtheorem{theoreme}{Th\'eor\`eme}[section]

\newtheorem{definition}[theoreme]{D\'efinition\rm}

\renewcommand{\theequation}{\arabic{equation}}
\def\og{\leavevmode\raise.3ex\hbox{$\scriptscriptstyle\langle\!\langle$~}}
\def\fg{\leavevmode\raise.3ex\hbox{~$\!\scriptscriptstyle\,\rangle\!\rangle$}}

\journal{the Acad\'emie des sciences}
\begin{document}
\centerline{Partial differential equations}
\begin{frontmatter}


\selectlanguage{english}
\title{Intrinsic Taylor formula for non-homogeneous Kolmogorov-type Lie groups
}


\selectlanguage{english}
\author[authorlabel1]{Stefano Pagliarani},
\ead{stepagliara1@gmail.com}
\author[authorlabel2]{Michele Pignotti}
\ead{michele.pignotti2@unibo.it}

\address[authorlabel1]{DEAMS, Universit\`a di Trieste, Via Tigor 22, 34124 Trieste, Italia. 
}
\address[authorlabel2]{Dipartimento di Matematica, Universit\`{a}
di Bologna, Piazza di Porta S. Donato 5, 40126 Bologna (Italy)}


\medskip
\begin{center}
{\small Received *****; accepted after revision +++++\\
Presented by £££££}
\end{center}

\begin{abstract}
\selectlanguage{english}
We prove an intrinsic Taylor-like formula for 
a class of 
Lie groups 
arising in the study of some sub-elliptic differential operators, namely the Kolmogorov operators. The estimate of the remainder is in terms of the intrinsic norm induced by such operators. These results extend the recent developments in \cite{PPP2016}, where a full characterization of the intrinsic H\"older spaces and their Taylor polynomials were given under the additional assumption that the Lie group is homogeneous in the sense of \cite{FollandStein1982}. Remarkably, the intrinsic Taylor polynomial admits the same representation as in the homogeneous case.



\vskip 0.5\baselineskip

\selectlanguage{francais}

\noindent{\bf R\'esum\'e} \vskip 0.5\baselineskip  

Nous d\'emontrons une formule intrins\`eque de type Taylor pour une classe de groupes de Lie d\'ecoulant de l'\`etude de certains op\'erateurs diff\'erentiels sous-elliptiques, \`a savoir les op\'erateurs de Kolmogorov. L'estimation du reste correspond \`a la norme intrins\`eque induite par ces op\'erateurs. Ces r\'esultats \'etendent les d\'eveloppements r\'ecents dans \cite{PPP2016}, o\`u une caract\'erisation compl\`ete des espaces de H\"older intrins\`eques et leurs polyn\^omes de Taylor ont \'et\'e donn\'es sous l'hypoth\`ese suppl\`ementaire que le groupe de Lie est homog\`ene dans le sens de \cite{FollandStein1982}. Remarquablement, le polyn\^ome intrins\`eque de Taylor admet la m\^eme repr\'esentation que dans le cas homog\`ene.

%
%
%


\end{abstract}
\end{frontmatter}


\selectlanguage{english}

\vspace{-10pt}
\section{Introduction}


Consider 
the non-commutative group $\Gc_B=(\Rdd,\circ
)$ defined by 
\begin{equation}\label{eq:translation}
 (t,x)\circ (s,\xi) = \left(t+s,{e^{s B}}x+\xi\right),\quad (t,x)^{-1}=\left(-t,{-e^{-tB}}x\right), \quad \text{Id}= (0,0),\quad (t,x),(s,\xi)\in \Rdd,
\end{equation}
where $B$ is a 
$(d\times d)$-square matrix such that
 \begin{equation}\label{eq:B_blocks_b}
\footnotesize
B=\left(
\begin{array}{ccccc}
B_{0,0}&B_{0,1}&\cdots&B_{0,r-1}&B_{0,r} \\ B_{1,0} & B_{1,1} &\cdots& B_{1,r-1} & B_{1,r} \\ 0 & B_{2,1} &\cdots& B_{2,r-1}& B_{2, r} \\
\vdots & \vdots &\ddots& \vdots&\vdots \\ 0 & 0 &\cdots& B_{r,r-1}& B_{r,r}
\end{array}
\right),\qquad 
\begin{array}{l}
\centerdot\ B_{i,j}\in \mathcal{M}^{p_i\times p_j}\\
\centerdot\ B_{j,j-1}\text{ have rank } p_j\\
\centerdot\ p_0\geq p_1\geq \cdots \geq p_r\geq 1\\
\centerdot\ p_0 + p_1 + \cdots + p_r = d.
\end{array}
\end{equation}
Hereafter, $\mathcal{M}^{p\times q}$ denotes the spaces of the $(p\times q)$-matrices with real entries, while $I_p$ will denote the $(p\times p)$ identity matrix. 
The group $\Gc_B$ plays a crucial role in the study of degenerate Kolmogorov  
operators
\begin{equation}\label{e1_b}
\Lc=\frac{1}{2}\sum_{i=1}^{p_{0}} A_{ij}(t,x) \p_{x_{i}x_{j}}+\langle Bx,\nabla_x\rangle+\p_{t},\qquad
  (t,x)\in\R\times\R^{d}, \qquad 1\le p_{0}\le d,
\end{equation}
with $A(t,x)\in\mathcal{M}^{p_0\times p_0}$, 
which is closely related to averaged-diffusion stochastic processes and whose applications include Mathematical Finance and Physics among others.  
In particular, $\Gc_B$ was first introduced in \cite{LanconelliPolidoro1994} as a group of left-invariant translations for the operator $\Lc_0$ obtained by $\Lc$ setting $A=I_{p_0}$. Defining
\begin{equation}\label{e3}
  X_{j}:=\p_{x_{j}},\quad j=1,\dots,p_{0},\quad\text{ and }\quad Y:=\langle
  Bx,{ \nabla_x}\rangle+\p_{t},
\end{equation}
in \cite{LanconelliPolidoro1994} it was proved that \eqref{eq:B_blocks_b} is equivalent to 
the H\"ormander's condition,
namely $\text{rank}\left(\text{Lie}(X_{1},\dots,X_{p_{0}},Y)\right)=d+1,$
which in turn implies that $L_0$ is hypoelliptic. Later (see \cite{Polidoro1994} and \cite{amrx}), hypoellipticity for a general $\Lc$ with $B$ as in \eqref{eq:B_blocks_b}, along with some regularity properties and Gaussian upper bounds for its fundamental solution, were proved under the assumption that $A(t,x)$ satisfies a uniform ellipticity condition on $\R^{p_0}$. In the above references, the classical notions of regularity based on the Euclidean distance are replaced with intrinsic notions of regularity related to the 
$B$-
semi-distance 
\begin{equation}\label{eq:semi_distance}
\|\z^{-1}\circ z\|_{B},\qquad z=(t,x),\z=(s,\xi)\in\Rdd,
\end{equation}
which was first introduced in \cite{LanconelliPolidoro1994}, where the $B$-norm $\|\cdot\|_B$ is defined as 
\begin{equation}\label{e7}
 \|{(t,x)}\|_B:=|t|^{1/2}+|x|_B,\qquad |x|_B:=
 \sum_{j=0}^r \sum_{i=\bar{p}_{j-1}+1}^{\bar{p}_j} |x_{i}|^{2j+1}, \qquad 
{\bar{p}_j:= \sum_{k=0}^j p_k} .
\end{equation}

It is now standard practice to conduct the study of $\Lc$-like operators 
by taking into account the intrinsic geometry induced by the operator, in particular by adopting suitable intrinsic H\"older spaces related to the semi-distance \eqref{eq:semi_distance}. Although  customized versions of such spaces were employed in several works (see 
\cite{Francesco}, 
\cite{Manfredini}, \cite{Lunardi1997}, 
\cite{NyPasPol10}, \cite{Priola} and \cite{Menozzi} among others), a complete characterization of the intrinsic H\"older spaces 
at any order and a systematic study of the related intrinsic Taylor polynomials was first performed in \cite{PPP2016}, 
under the additional assumption that the blocks $(B_{i,j})_{i\leq j}$ in
\eqref{eq:B_blocks_b} are null. {The latter} study is crucial in order to derive high-order results for the solutions of $\Lc$, e.g. Schauder estimates
, or asymptotic expansions (see \cite{PagliaraniPascucciPignotti2017}).  
In 
\cite{PPP2016}, 
the theory of intrinsic regularity was developed so as to mimic the general classical theory, namely: first one defines, recursively, the intrinsic H\"older space $C^{n,\alpha}_B$ by only specifying the regularity along the vector fields $X_1,\cdots,X_{p_0},Y$, and then one proves an $n$-th order Taylor formula for a function $u\in C^{n,\alpha}_B$ with remainder expressed in terms of the intrinsic semi-distance \eqref{eq:semi_distance}. Here we extend the program successfully pursued in \cite{PPP2016} to the case of a general matrix $B$ satisfying \eqref{eq:B_blocks_b}; i.e.:
\begin{enumerate}
\item give a new definition of intrinsic $C_B^{n,\a}$ regularity for $n\in \N_0$, $\a\in ]0,1]$; this appears to be the minimal regularity required in order
for a Taylor formula to hold with an estimate of the remainder proportional to $\|u\|_{C^{n,\a}_{B}}\|\z^{-1}\circ z\|_{B}^{n+\a}$;
\item prove a novel and explicit expression for the Taylor polynomial 
$T_n u(\z,\cdot)$ 
 centered in $\z=(s,\x)\in \Rdd $ of a function $u\in C_B^{n,\a}$,  namely
\begin{equation}\label{eq:ste_Tay_pol}
 T_n u(\z,z):=  \sum_{\substack{k\in\N_0,\, \beta\in\N_0^d\\0\leq 2 k + |\beta|_B \leq n}}\frac{1}{k!\,\beta!}
 \big( Y^k \partial_{\xi}^{\beta}u(s,\xi)\big) (t-s)^k\big( x-e^{(t-s)B}\xi  \big)^{\beta},\qquad
 z=(t,x)\in\Rdd, 
\end{equation}
{where 
$|\beta|_B:  =   \sum_{j=0}^r \sum_{i=\bar{p}_{j-1}+1}^{\bar{p}_j} (2 j +1) \beta_i$ is the $B$-length of $\beta$, 
$\partial_{\xi}^{\beta}:=\partial_{\xi_1}^{\beta_1} \partial_{\xi_2}^{\beta_2}\cdots\partial_{\beta_d}^{\xi_d}$ is a Euclidean multi-derivative, whereas $Y^{k}$ is meant as a $k$-th order Lie derivative.}
\end{enumerate}

We recall that, from the algebraic stand-point, the hypothesis $B_{i,j}=0$, $i\leq j$, previously employed in \cite{PPP2016} is equivalent to saying that $\Gc_B$ is homogeneous (in the sense of \cite{FollandStein1982}) with respect to the family of automorphisms $\left(D({\lambda})\right)_{\lambda>0}$, called dilations, defined as $
 D(\lambda)=\textrm{diag}\big(\lambda^2,\lambda I_{p_0},\lambda^{3}I_{p_1},\cdots,\lambda^{2r+1}I_{p_r}\big)$, 
with respect to which the norm $\|\cdot\|_B$ is homogeneous of degree one. Due to this reason the latter is sometimes referred to as \emph{homogeneous} norm.  
Previous results about intrinsic Taylor polynomials on homogeneous Lie groups were proved in greater generality in \cite{FollandStein1982}, and also recently in \cite{Bonfiglioli2009}. 
However, in these references the Taylor polynomials were defined for functions that are differentiable in the Euclidean sense; consequently the constants appearing in the estimates of the remainder would depend on  the Euclidean H\"older norms. {Moreover, the number of terms appearing in the Taylor polynomials grew exponentially with $n$, while in \eqref{e7} it grows only linearly.} 
Even more importantly, those results strongly rely on the use of the dilations and thus on the homogeneity assumption on $\Gc_B$. 
By opposite, the proof technique adopted in \cite{PPP2016} owes to the fact that the $B$-intrinsic norm in \eqref{e7} is somehow well-behaved with respect to the \emph{stratification} on 
\begin{align}
\Rdd &= \Rpp \times \cdots \times \R^{p_r}  \\
&\cong \underbrace{\text{span}(Y)|_{x=0}}_{=:W_{\text{time}}} \oplus \underbrace{\text{span}(X_1,\cdots,X_{p_0})}_{=:W_0} \oplus    \underbrace{[W_0,W_{\text{time}}]|_{x=0}}_{=:W_1} \oplus \cdots \oplus   \underbrace{[W_{r-1},W_{\text{time}}]|_{x=0}}_{=:W_{r}} 
\end{align}
equipped with the \emph{formal degrees} $\mathfrak{d}(W_0)$ and $\mathfrak{d}(W_{\text{time}})$ equal to $1$ and $2$ respectively, and the formal degrees $\mathfrak{d}(W_{i})=\mathfrak{d}(W_{i-1})+\mathfrak{d}(W_{\text{time}})$ for $i=1,\cdots,r$; the extension to the non-homogeneous case relies on the fact that $\|\cdot \|_B$ remains well-behaved w.r.t. to this stratification. 
Finally note that, remarkably, $T_n u(\z,\cdot)$ admits the same representation as in the homogeneous case, though it ceases to be a polynomial w.r.t. the time-variable as the matrix $B$ is in general not nilpotent.

\section{{
Intrinsic H\"older spaces and Taylor formula}}\label{sec2}

We start by introducing the notions of $B$-intrinsic regularity {in terms of suitable} $B$-intrinsic H\"older spaces. 
For any $z=(t,x)\in\Rdd$, we denote by 
\begin{equation}\label{eq:def_curva_integrale_campo}
 e^{\d X_{i} }(t,x)=(t,x+\delta e_i),\quad i=1,\cdots,p_0,\qquad
 e^{\d Y }(t,x)=(t+\delta,e^{\delta B}x),\qquad \delta>0,
\end{equation}
the integral curves of the fields $X_1,\cdots,X_{p_0},Y$ starting at $z$. {Here $e_i$ denotes the $i$-th element of the canonical basis of $\R^d$.}
 Now, let $\O$ be a domain in $\Rdd$. For any $z\in\O$ we set
  $$\d_{z}:=\sup\left\{\bar{\d}\in\,]0,1]\mid e^{\d X_1}(z),\cdots,e^{\d X_{p_0}}(z),e^{\d Y}(z)\in\O\text{ for any }\d\in [-\bar{\d},\bar{\d}]\right\}.$$
If $\O_{0}$ is a bounded domain with $\overline{\O}_{0}\subseteq\O$, we set
  $\d_{\O_{0}}=\min_{z\in \overline{\O}_{0}}\d_{z}.$
Note that $\d_{\O_{0}}\in\,]0,1]$. 
\begin{definition}\label{def:intrinsic_alpha_Holder2}
For $\a_1\in\,]0,1]$ and $\a_2\in]0,2]$, we say that $u\in C_{X_i}^{\alpha_1}(\O)$, $i=1,\cdots,p_0$, and $v\in C_{Y}^{\alpha_2}(\O)$, if for any
bounded domain $\O_{0}$ with $\overline{\O}_{0}\subseteq\O$, the following semi-norms are finite:
\begin{equation}\label{e12}
 \left\|u\right\|_{C_{X_i}^{\alpha_1}(\O_{0})}:=\!\!\sup_{z\in \O_{0}\atop 0<|\d|<\d_{\O_{0}}} \!\!\frac{\left|u\left(e^{\delta X_i }(z)\right)-
 u(z)\right|}{|\delta|^{\a_1}}, \quad i=1,\cdots, p_0, \qquad \left\|v\right\|_{C_{Y}^{\alpha_2}(\O_{0})}:=\!\!\sup_{z\in \O_{0}\atop 0<|\d|<\d_{\O_{0}}}\!\! \frac{\left|v\left(e^{\delta Y }(z)\right)-
 v(z)\right|}{|\delta|^{\frac{\a_2}{2}}}.
\end{equation}
\end{definition}
We can now define the $B$-intrinsic H\"older spaces on 
$\mathcal{G}_B$ ($B$-H\"older spaces).
\begin{definition}\label{def:C_alpha_spaces}
Let $\a\in\,]0,1]$, then:
\begin{itemize}
  \item [i)] $u\in C^{0,\a}_{B}(\O)$ if $u\in C^{\a}_{Y}(\O)$ and $u\in C^{\a}_{{X_i}}(\O)$ for any $i=1,\dots,p_{0}$; 
  \item [ii)] $u\in C^{1,\a}_{B}(\O)$ if $u\in C^{1+\a}_{Y}(\O)$ and $\p_{x_{i}}u\in C^{0,\a}_{B}(\O)$ for any
  $i=1,\dots,p_{0}$; 
  \item [iii)] {For} $k\in\Nb$ with $k\ge2$, $u\in C^{k,\a}_{B}(\O)$ if $Yu\in C^{k-2,\a}_{B}(\O)$ and $\p_{x_{i}}u\in C^{k-1,\a}_{B}(\O)$ for any
  $i=1,\dots,p_{0}$. 
\end{itemize}
Moreover, the space $C^{k,\a}_{B}(\O)$ is equipped, for any bounded domain ${\O}_{0}\subset\O$, with the $\O_0$-seminorm
\begin{equation}
\|{u}\|_{C^{k,\a}_{B}(\O_0)}:=\left\{
\begin{aligned}
&\|{u}\|_{C^{\a}_{Y}(\O_0)}+\sum_{i=1}^{p_0} \|{u}\|_{C^{\a}_{X_i}(\O_0)},&& k=0\\
&\|{u}\|_{C^{\a+1}_{Y}(\O_0)}+\sum_{i=1}^{p_0} \|{\partial_{x_i}u}\|_{C^{0,\a}_{B}(\O_0)}, && k=1\\
&\|{Y u}\|_{C^{k-2,\a}_{B}(\O_0)}+\sum_{i=1}^{p_0} \|{\partial_{x_i}u}\|_{C^{k-1,\a}_{B}(\O_0)}, && k\geq 2.
\end{aligned}
\right.
\end{equation}

\end{definition}
{Recall that, in the previous definition as well as in the whole paper, $Y u$ is meant as a Lie derivative.} 
\begin{theorem}\label{th:main}
Let $\O$ be a domain of $\Rdd$, $\alpha\in ]0,1]$ and $n\in\Nzero$. If $u\in C^{n,\a}_{B}(\O)$ then it holds:
\begin{enumerate}
\item[1)] there exist 
\begin{equation}\label{eq:maintheorem_part1}
 Y^k \partial_x^{\beta}u\in C^{n-2k-|\beta|_B,\alpha}_{B}(\O),
 \qquad {0\leq 2 k + |\beta|_B \leq n};
\end{equation}
\item[2)]
for any $\z_0\in\Omega$, there exist two bounded domains $U,V$, such that $\zeta_0\in U\subset V\subset\O$ 
and
\begin{equation}\label{eq:estim_tay_n_loc}
 \left|u(z)-T_n u(\z,z)\right|\le c_{B,U
 } \|u\|_{C^{n,\a}_{B}( V
 )}\|\z^{-1}\circ z\|_{B}^{n+\a},\qquad z,\zeta \in U,
\end{equation}
where $c_{B,U}$ is a positive constant 
and $T_n u(\z,\cdot)$ is the
\emph{$n$-th order $B$-intrinsic Taylor polynomial of $u$ centered in $\z
$} as defined in \eqref{eq:ste_Tay_pol}.
\end{enumerate}
\end{theorem}

\section{Proof of Theorem \ref{th:main}}\label{sec3}
For sake of brevity, we only prove the statement for $r=1$, which is $B=
(B_{i,j})_{i,j\in\{0,1\}}$ with $B_{i,j}\in \mathcal{M}^{p_i\times p_j}$ and $B_{1,0}$ has full rank.  This case is complex enough to see the conceptual difficulties that arise from dropping the homogeneity assumption on $\mathcal{G}_B$. On the other hand, the proof for a general $r\geq 1$ is only a lengthy and technical extension.
\begin{notation}
Throughout this section we will use the notation $z=(t,x,y)$ or $\zeta=(s,\xi,\eta)$ to indicate a general element of $\Rpp$. 
 Moreover, we will denote by $c$ any positive constant that depends on $B$ and on the domain $U$ in Theorem \ref{th:main}, at most. 
\end{notation}

The first task is connecting two points in $\Rpp$ using integral curves. 
To obtain an increment in the $x$-variables it is enough to move along the integral curves of the fields $X_1,\cdots,X_{p_0}$, i.e.
\begin{equation}
\g^{(0)}_{v,\d}(t,x,y) : = 
(t, x + \d v, y), \qquad 
v\in\Rpz, \ 
{\delta\in\R}.
\end{equation}
To understand how to obtain an increment in the $y$-variables, it is useful to observe that 
\begin{equation}
\label{eq:par_y}
[v_1 X_1 + \cdots + v_{p_0} X_{p_0}, Y] - \langle {\nabla_x}, B_{0,0} v \rangle = \langle {\nabla_y}, B_{1,0} v \rangle,\qquad v\in\Rpz.
\end{equation}
It is thus reasonable to approximate the integral curves of the vector field on the right-hand side as
\begin{align}
\g_{v,\d}(t,x,y)& : =
\g^{(0)}_{B_{0,0} v,-\d^3}\left(e^{-\d^{2}Y}\left(\g^{(0)}_{v,-\d}\left(e^{\d^{2}Y}\left(\g^{(0)}_{v,\d}(t,x,y)\right)\right)\right)\right)
\\ 
&  =\big( t, x ,  y+\d^3 B_{1,0} v   \big)  - \d^5 \bigg(0 , \sum_{n=0}^{\infty} \frac{(-1)^{n}\d^{2n}}{(n+2)!} B^{n+2} (v,0)^\top \bigg), \qquad v\in\Rpz, \ 
{\delta\in\R}.
\label{eq:def_gamma1}\end{align}
The leading order increment is proportional to $\d^3 
$, along the $y$ variable only. However,
due to the non-homogeneous structure of $B$ (the block $B_{0,0}$ is not null), the higher order increment affects both the components $x$ and $y$.
To correct this, we employ again the curve $\g^{(0)}$. Set
\begin{align}
\label{eq:def_g1}
g_{v,\d}(t,x,y)& : = \g^{(0)}_{v',\d'}\left(\g_{v,\d}(t,x,y)\right), \qquad v\in\Rpz, \ 
{\delta\in\R},
\end{align}
where 
\begin{equation}
\label{eq:def_g2}
v' = v'(\d,v) = \sum_{n=0}^{\infty} \frac{(-1)^{n}\d^{2n}}{(n+2)!} B^{n+2}_{0,0} v , \quad \d' =\d'(\d) =  \d^5,
\end{equation}
and $B^{n+2}_{0,0}$ is the top-left $(p_0\times p_0)$-submatrix of $B^{n+2}$.
\begin{lemma}\label{lem:connect}
There exists $\eps>0$, only dependent on $B$, such that: for any $\eta\in \Rpo$ with $|\eta|\leq \eps$, there exist $v\in \Rpz$ with $|v|=1$ and $\delta\geq 0$ such that 
\begin{equation}
\label{eq:g1_estim}
g_{v,\d}(t,x,y) = (t, x, y + \eta ), \qquad \text{and} \qquad |\d| \leq c
 |\eta|
 ^{\frac{1}{3}}.
\end{equation}
\end{lemma}
\noindent{\bf Proof.}
By 
{\eqref{eq:def_g1} and \eqref{eq:def_gamma1}} we obtain
\begin{align}
g_{v,\d}(t,x,y) - (t, x, y ) 
= \big( 0, 0, \d^3  R(\delta,v) \big), && R(\delta,v):= \sum_{n=0}^{\infty} \frac{(-1)^{n}\d^{2n}}{(n+1)!} B^{n+1}_{1,0} v,
\end{align}
where $B^{n+1}_{1,0}$ is the bottom-left $(p_1\times p_0)$-submatrix of $B^{n+1}$. Therefore, denoting by $\mathbb{S}^{p_0-1}$ the unitary sphere in $\Rpz$, we have to find some $(\d,v)\in[0,\infty[\times\mathbb{S}^{p_0-1}$ that solves the equation
\begin{equation}\label{eq:algebraic_eq}
\d^3 R(\delta,v) = \eta.
\end{equation} 
Since $B_{1,0}$ has full rank it is not restrictive to assume $p_0=p_1$, and thus $B_{1,0}$ invertible. In particular, $R(0,v)=B_{1,0}v$, which implies that $R(0,\cdot)$ is a bijective and linear function. Moreover, since $R$ is globally $C^1$, there exists $\bar{\d}>0$ such that $R(\d,\cdot)$ is still a bijective linear function for any $\delta\leq \bar{\d}$. 
In particular, when restricted to $\mathbb{S}^{p_0-1}$, $(\d^3 R(\d,\cdot))_{0\leq\d\leq\bar{\d}}$ is a continuous family of embeddings that collapses to zero at $\d=0$. Therefore, by employing Jordan-Brouwer's Theorem, it is possible to prove that Eq. \eqref{eq:algebraic_eq} admits a solution $\big(\d(\eta),v(\eta)\big)\in[0,\bar{\d}]\times \mathbb{S}^{p_0-1}$ for any $|\eta|\leq \eps$, where $\eps>0$ only depends on $B$. We now prove the second part of \eqref{eq:g1_estim}. Choosing $\eps$ small enough, it holds $\big|R\big(\delta(\eta),v(\eta)\big)\big|\in [\norm{B_{1,0}}-\eps, \norm{B_{1,0}}+\eps] $, and 
by \eqref{eq:algebraic_eq}, 
\begin{equation}
\big|\d(\eta)\big|^3  = \frac{ |\eta| }{ \big| R\big(\delta(\eta),v(\eta)\big) \big|} \leq \frac{ |\eta| }{\max(0,\norm{B_{1,0}}-\eps)}.
\end{equation} 
Again, taking $\eps$ suitably small yields the result. 
\hfill$\Box$

We are now in the position to prove Theorem \ref{th:main}. 

\noindent\textbf{Proof of Theorem \ref{th:main}
}.
Analogously to the homogeneous setting, the cases $n=0,1,2,3$ have to be proved separately, while for $n>3$ the proof is by induction on $n$. For sake of brevity, here we only provide a proof for $n=0$ and $n=3$, these being the most interesting and difficult steps. On the one hand, the proof for $n=0$ allows to appreciate how the connection Lemma \ref{lem:connect} along with the regularity along the fields can be used, in a rather simple way, in order to obtain the most basic result, namely the H\"olderianity with respect to the $B$-intrinsic norm. On the other hand, the proof for $n=3$ enlightens the main difficulty of the whole proof, namely proving the existence of the first order partial derivative w.r.t. the $y$-variable. Note that the existence of the latter is not trivially ensured by the definition of $C^{3,\a}_{B}(\O)$, as the existence of $X_i Y u$ and $Y X_i$, and thus the commutators $[X_i,Y]$, are only meant in the sense of Lie derivatives. As for the steps $n=1$, $n=2$, these are just simplifications of the case $n=3$, whereas the inductive step for $n>3$ is totally analogous to the homogeneous case. 


\smallskip\noindent\underline{\emph{Case $n=0$}}:
{We only need to prove Part 2).} Let $U\subset \Omega$ be a bounded domain suitably small so as to ensure that all the integral curves that are employed below to connect $z$ and $\z$ are entirely contained in the bounded domain $V\subset \Omega$. 
The first step is to {bound the increment w.r.t. the time variable by employing the integral curve of $Y$ in \eqref{eq:def_curva_integrale_campo}:}
\begin{align}\label{ande24}
  \left|u(t,x,y)-u(s,\x,\y)\right|&\le    \big| u(t,x,y)-{u\big( e^{(t-s)Y}(s,\x,\y) \big)}\big|+\big|{u\big( e^{(t-s)Y}(s,\x,\y)\big)}-u(s,\x,\y)\big|\\
  &\le \big|u(t,x,y)-u\big(t,e^{(t-s)B}(\x,\y)^\top\big)\big| + c
\, \norm{u} _{C^{\a}_Y(V)} |s-t|^{\frac{\a}{2}},
\end{align}
where we used triangular inequality in the first line and $u\in C^{\a}_Y(\Omega)$ in the second.
Note that $\big|(x,y)^\top-e^{(t-s)B}(\x,\y)^\top\big|_{B}^{\a}\le \|\z^{-1}\circ z\|_B^{\a}$ and 
thus we only need to prove
\begin{equation}\label{eq:ste23}
 |u(t,\x,\y)-u(t,x,y)|\leq c\, \norm{u}_{C_B^{0,\a}(V)}(|\x-x|+|\y-y|^{\frac{1}{3}})^{\a}.
 \end{equation} 
We can use again triangular inequality and write 
\begin{align}
u(t,\x,\y)-u(t,x,y) &= \big(u(t,\x,\y)-u(t,x,\y)\big) +\big(u(t,x,\y)-u(t,x,y)\big)\\
& = \big(u(t,\x,\y)-u(t,x,\y)\big) +\Big(u\big(g_{\d,v}(t,x,y)\big)-u(t,x,y)\Big)\quad \text{(by Lemma \ref{lem:connect})}
\end{align}
with $|v|=1$ and ${|\d|\leq c |\eta-y|^{\frac{1}{3}}}$. By using $u\in C^{\a}_{{X_i}}(\Omega)$, $i=1,\dots,p_0$, in order to bound the first term, 
together with $u\in C^{\a}_{{Y}}(\Omega)$ to bound the second, we obtain \eqref{eq:ste23}, which concludes the proof for $n=0$.

\smallskip\noindent\underline{\emph{Case $n=3$}}: To shorten notation we only prove the case $p_0=p_1=1$. The difficulty in considering multi-dimensional blocks is purely notational. We first prove Part 1). Fix an arbitrary bounded domain $\Omega_0\subset\Omega$. Proceeding 
as in 
the homogeneous case (see the proof of \cite[Proposition 5.25]{PPP2016}), one obtains
\begin{align}
\label{eq:estim_gamma}
\big| u\big( \gamma^{(i)}_{\d,v}(z) \big) - \bar{T}_3 u\big(z, \gamma^{(i)}_{\d,v}(z) \big)  \big| & \leq c\, \|u\|_{C^{3,\a}_{B}(\Omega_0)} |\d|^{3+\a},\qquad i=0,1,
\end{align}
for any $z\in \Omega_0$, and $v,\d\in\R$ with $|v|=1$ and $|\delta|$ suitably small, where we set
\begin{equation}
\bar{T}_3 u(z,\zeta) = \sum_{i=0}^3 \frac{ (\x-x)^i }{i!} \partial^i_x u(z)+\frac{\eta - y}{B_{1,0}} \big( [\p_{x}, Y] -B_{0,0}\p_x\big) u (z) , \qquad z=(t,x,y),\,\zeta=(t,\x,\eta).
\end{equation}
The last term in the right-hand side is inspired by \eqref{eq:par_y} to mimic a {partial derivative w.r.t. y} 
and is well defined when applied to $u\in C^{3,\a}_B(\Omega)$.
We now prove 
\begin{align}
\label{eq:estim_g}
\big| u\big( g_{\d,v}(z) \big) - \bar{T}_3 u\big(z, g_{\d,v}(z) \big)  \big| & \leq c\, \|u\|_{C^{3,\a}_{B}(\Omega_0)} |\d|^{3+\a},
\end{align}
where $g_{v,\d}$ is as defined in \eqref{eq:def_g1}-\eqref{eq:def_g2}. 
Setting $z':=\g_{\d,v}(z)$ and $z''=g_{v,\d}(z)$ we have
\begin{align}
&u\big( z'' \big) - \bar{T}_3 u\big(z, z'' \big)= F_1 +F_2,\\
F_1&=\left(u( z'') - \bar{T}_3 u(z', z'' )\right) + \left(u(z') - \bar{T}_3 u(z, z')\right),\\
F_2&=   \bar{T}_3 u(z', z'')+\bar{T}_3 u(z,z') -u(z')-\bar{T}_3 u(z, z'') .
\end{align}
Now, \eqref{eq:estim_gamma} and \eqref{eq:def_g2} yield $|F_1|\leq c\, \|u\|_{C^{3,\a}_{B}(\Omega_0)} |\d|^{3+\a}$; as for $F_2$ it holds:
 \begin{equation}
 |F_2|=\left|\sum_{i=1}^3 \frac{1}{i!} \left(\p_{x}^i u(z')-\p_{x}^i u(z)\right)(\d' v')^i\right|\leq c\, \|u\|_{C^{3,\a}_{B}(\Omega_0)} |\d' v'|\leq c\, \|u\|_{C^{3,\a}_{B}(\Omega_0)}|\d|^{3+\a},
 \end{equation}
where we used {$\partial^i_x u \in C^{3-i,\alpha}_B(\Omega)$ and Theorem \ref{th:main} for $n=0,1,2,$} 
to prove the first inequality, and \eqref{eq:def_g2} to prove the second one. This proves \eqref{eq:estim_g}. 
We are now able to prove 
differentiability along the $y$ direction. For any $z=(t,x,y)\in\Omega_0$ and $\eta\in\R$ with $|\eta|$ small enough, choosing $v$ and $\d$ as given by Lemma \ref{lem:connect} yields
\begin{align}
\label{eq:estimate_gamma2}
\big|u(t,x,y+\eta)-\bar{T}_3 u\big((t,x,y),(t,x,y+\eta)\big)\big|&=\big|u\big(g_{v,\d}(t,x,y)\big)- \bar{T}_3 u\big((t,x,y), g_{\d,v}(t,x,y) \big)\big|\\
&\leq c \|u\|_{C^{3,\a}_{B}(\Omega_0)} |\d|^{3+\a}\leq  c \|u\|_{C^{3,\a}_{B}(\Omega_0)} |\eta|^{1+\frac{\a}{3}},
\end{align}
where we used \eqref{eq:estim_g} in to obtain the first inequality, and \eqref{eq:g1_estim} to obtain the second.
Thus $\p_y u(z)$ exists and 
\begin{equation}
\p_y u(t,x,y)= \frac{1}{B_{1,0}} \big( [\p_{x}, Y] -B_{0,0}\p_x\big) u (t,x,y).
\end{equation}
Furthermore, $u\in C^{3,\alpha}_B(\Omega)$ implies $\p_y u\in C^{0,\alpha}_B(\Omega)$, which is Part 1) of Theorem \ref{th:main} for $n=3$.

The proof of Part 2) is analogous 
to the homogeneous case treated in \cite[Subsections 5.2.2, 5.2.3]{PPP2016}. For brevity, we 
only give an account of the main idea. The first step is to prove 
\eqref{eq:estim_tay_n_loc} for $z,\z$ sharing the same time component, say $z=(t,x,y)$, $\z=(t,\x,\y)$. We define the intermediate point $z_1:=(t,\x,y)$ and write 
\begin{equation}
\label{eq:esti_tay1}
u(z)-T_3 u (\z,z)= \left(u(z)-T_3 u (z_1,z)\right)+\left(T_3 u (z_1,z)-T_3 u (\z,z)\right)
.
\end{equation}
Now, an application of the mean value theorem along $x$  in 
the first addend, and along $y$ in 
in the second, yields the desired estimate. 
Note that this step strongly relies on Part 1), namely the differentiability w.r.t. the $y$ variable.
Finally, the 
case of general $z,\zeta$ can be {reduced to the above one} similarly to \eqref{ande24}. 

\hfill$\Box$
\begin{footnotesize}
\bibliography{Bibtex-Master-3.00}
\end{footnotesize}
\end{document}